\documentclass[a4paper,11pt,leqno]{article}%
\usepackage{graphicx}
\usepackage{amsmath}
\usepackage{amsfonts}
\usepackage{amssymb}%
\newtheorem{theorem}{Theorem}

\newtheorem{propdef}[theorem]{Proposition-Definition}

\newtheorem{notation}{Notation}

\newtheorem{proposition}[theorem]{Proposition}
\newtheorem{remark}{Remark}

\newenvironment{proof}[1][Proof]{\textbf{#1.} }{\ \rule{0.5em}{0.5em}}
\begin{document}

\title{Some properties of $(\mathcal{C},\mathcal{E},\mathcal{P})$-algebras :
Overgneration and $0$-order estimates.}
\author{Antoine Delcroix \thanks{Equipe AANL, Laboratoire AOC -- Campus de Fouillole,
Universit\'{e} des Antilles et de la Guyane, 97159 Pointe-\`{a}-Pitre,
Guadeloupe (France), E-mail: antoine.delcroix@univ-ag.fr}}
\maketitle

\begin{abstract}
We give a new definition of the so-called overgenerated rings, which are the
usual tool used to define the asymptotic structure of a $(\mathcal{C}%
,\mathcal{E},\mathcal{P})$-algebra, written as a factor space $\mathcal{M}%
_{(A,\mathcal{E},\mathcal{P})}/\mathcal{N}_{(I_{A},\mathcal{E},\mathcal{P})}$.
With this new definition and in the particular case of $\mathcal{E}%
=\mathrm{C}^{\infty}$, we show that a moderate element i.e. in $\mathcal{M}%
_{(A,\mathcal{E},\mathcal{P})}$ is negligible if and only if it satisfies the
$\mathrm{C}^{0}$-order estimate for the ideal $\mathcal{N}_{(I_{A}%
,\mathcal{E},\mathcal{P})}$.\ 

\end{abstract}

\noindent\textbf{Mathematics Subject Classification (2000):} 35A20, 35A25,
35D05, 46F30, 46T30\textbf{\smallskip}


\section{Introduction}

More than ten years ago, J.-A. Marti has introduced the structure of
$(\mathcal{C},\mathcal{E},\mathcal{P})$-algebras \cite{JAM1} which is a
refinement of the Colombeau simplified algebras of new generalized functions
\cite{Col1,GKOS,Ober1}. These $(\mathcal{C},\mathcal{E},\mathcal{P})$-algebras
are constructed on a base (pre)sheaf $\mathcal{E}$, which is usually a
(pre)sheaf of algebras equipped with a topology $\mathcal{P}$ and on an
asymptotic structure given by a ring $\mathcal{C}$ of generalized constants.
Except in a few cases (for example \cite{BENU1}), the sheaf $\mathcal{E}$ is
chosen to be the sheaf of smooth functions. Roughly speaking a presheaf of
$(\mathcal{C},\mathcal{E},\mathcal{P})$-algebras is a presheaf of factor
algebras $\mathcal{M}_{(A,\mathcal{E},\mathcal{P})}/\mathcal{N}_{(I_{A}%
,\mathcal{E},\mathcal{P})}$, where $\mathcal{M}_{(A,\mathcal{E},\mathcal{P})}$
(resp $\mathcal{N}_{(I_{A},\mathcal{E},\mathcal{P})}$) is the (pre)sheaf of
algebras of moderate elements (resp the (pre)sheaf of ideals of negligible
elements) over the factor ring $\mathcal{C}=A/I_{A}$.\ (The moderateness and
the negligibility are defined by the asymptotic structure given by the ring
$\mathcal{C}$.)\smallskip

Although these $(\mathcal{C},\mathcal{E},\mathcal{P})$-type algebras have
proved their efficiency to give a meaning and to solve singular differential
problems, the investigations on their intrinsic properties have not yet been
developed, expect first attempts concerning topology in \cite{ADDS2} and
asymptotic analysis in \cite{DEMAROB}. In this paper, after recalling for sake
of self contentedness the basic notions on $(\mathcal{C},\mathcal{E}%
,\mathcal{P})$-type algebras (Section \ref{CEP08SecCEP}), we present a new
definition of the so-called overgenerated rings in Section
\ref{CEP08SecOverTop}.\ This concept of overgeneration is one essential part
of the theory for it allows to adapt the algebraic structure to the
singularities of the problems. Taking advantage of this new definition and
analogously to the theorem 1.2.3. of \cite{GKOS} for Colombeau simplified
algebras, we show in Section \ref{SecAL} that a moderate element (i.e. in
$\mathcal{M}_{(A,\mathcal{E},\mathcal{P})}$) is negligible (i.e. in
$\mathcal{N}_{(I_{A},\mathcal{E},\mathcal{P})}$) if and only if it satisfies
the $\mathrm{C}^{0}$-order estimate of the ideal, for the case of
$\mathcal{E}=\mathrm{C}^{\infty}$ endowed with its classical topology. This
property allows radical simplification of many proofs of existence and
uniqueness for differential problems.

\section{The presheaf of $(\mathcal{C},\mathcal{E},\mathcal{P})$-type
algebras\label{CEP08SecCEP}}

%

%TCIMACRO{\TeXButton{Compteurequation}{\setcounter{equation}{6}}}%
%BeginExpansion
\setcounter{equation}{6}%
%EndExpansion
We begin by recalling the algebraical construction of $(\mathcal{C}%
,\mathcal{E},\mathcal{P})$-algebras \cite{JAM1,JAM2} as improved in
\cite{DEMAROB}. Let:

\noindent$\left(  1\right)  $~$\Lambda$ be a directed set with partial order
relation $\preceq$;

\noindent$\left(  2\right)  ~A$ be a solid subring of the ring $\mathbb{K}%
^{\Lambda}$\ ($\mathbb{K=R}$ or $\mathbb{C}$): whenever $(\left\vert
s_{\lambda}\right\vert )_{\lambda}\leq(\left\vert r_{\lambda}\right\vert
)_{\lambda}$\ for some $\left(  (s_{\lambda})_{\lambda},(r_{\lambda}%
)_{\lambda}\right)  \in\mathbb{K}^{\Lambda}\times A$, that is, $\left\vert
s_{\lambda}\right\vert \leq\left\vert r_{\lambda}\right\vert $ for all
$\lambda$, it follows that $(s_{\lambda})_{\lambda}\in A$;

\noindent$\left(  3\right)  $~$I_{A}$ be a solid ideal of $A$;

\noindent$\left(  4\right)  $~$\mathcal{E}$ be a sheaf of $\mathbb{K}%
$-topological algebras over a topological space $X$.\smallskip

Suppose that for any open set $\Omega$ in $X$, the topology of the algebra
$\mathcal{E}(\Omega)$ is defined by a family $\mathcal{P}(\Omega)$ of
seminorms such that:

\noindent$\left(  5\right)  $~Whenever $\Omega_{1}$, $\Omega_{2}$ are two open
subsets of $X$ with $\Omega_{2}\subset\Omega_{1}$ and $\rho_{2}^{1}$ is the
restriction operator $\mathcal{E}(\Omega_{1})\rightarrow\mathcal{E}(\Omega
_{2})$, then, for each $p_{2}\in\mathcal{P}(\Omega_{2})$, the seminorm
$p_{1}=p_{2}\circ\rho_{2}^{1}$ extends $p_{2}$ to $\mathcal{P}(\Omega_{1})$;

\noindent$\left(  6\right)  $~Whenever $\Theta=(\Omega_{h})_{h\in H}$ is a
family of open sets in $X$ with $\Omega=\cup_{h\in H}\Omega_{h}$, then, for
each $p\in\mathcal{P}(\Omega)$, there exist a finite subfamily $\left(
\Omega_{i}\right)  _{1\leq i\leq n}$ of $\Theta$ and corresponding seminorms
$p_{i}\in\mathcal{P}(\Omega_{i})$, $1\leq i\leq n$, such that%
\[
\forall u\in\mathcal{E}(\Omega),\ \ p\left(  u\right)  \leq p_{1}\left(
u\left\vert _{\Omega_{1}}\right.  \right)  +\ldots+p_{n}(u|_{\Omega_{n}}).
\]

Define $\mathcal{C}=A/I_{A}$ and $\left\vert B\right\vert =\left\{  \left(
\left\vert r_{\lambda}\right\vert \right)  _{\lambda},\ \left(  r_{\lambda
}\right)  _{\lambda}\in B\right\}  $ ($B=A$ or $I_{A}$).\ From $(2)$, it
follows that $\left\vert A\right\vert $ is a subset of $A$ and that
$A_{+}=\{(b_{\lambda})_{\lambda}\in A,\,\forall\lambda\in\Lambda,\ b_{\lambda
}\geq0\}=\left\vert A\right\vert $. The same holds for $I_{A}$. Furthermore,
$(2)$ implies also that $A$ is a $\mathbb{K}$-algebra \cite{DEMAROB}.\ With
these notations, set%
\begin{align*}
\mathcal{M}(\Omega)  &  =\mathcal{M}_{(A,\mathcal{E},\mathcal{P})}%
(\Omega)=\left\{  (u_{\lambda})_{\lambda}\in\left[  \mathcal{E}(\Omega
)\right]  ^{\Lambda}\mid\forall p\in\mathcal{P}(\Omega),\ \left(
(p(u_{\lambda})\right)  _{\lambda}\in\left\vert A\right\vert \right\}  ,\\
\mathcal{N}(\Omega)  &  =\mathcal{N}_{(I_{A},\mathcal{E},\mathcal{P})}%
(\Omega)=\left\{  (u_{\lambda})_{\lambda}\in\left[  \mathcal{E}(\Omega
)\right]  ^{\Lambda}\mid\forall p\in\mathcal{P}(\Omega),\ \left(
p(u_{\lambda})\right)  _{\lambda}\in\left\vert I_{A}\right\vert \right\}  .
\end{align*}

\begin{propdef}
\label{SPSPropSubAlg}\cite{DEMAROB, JAM1}~\newline$\left(  i\right)
$~$\mathcal{M}_{(A,\mathcal{E},\mathcal{P})}$\textit{\ (resp. }$\mathcal{N}%
_{(I_{A},\mathcal{E},\mathcal{P})}$) \textit{is a sheaf of }$\mathbb{K}%
$-\textit{subalgebras (resp. of ideals) of the sheaf }$\mathcal{E}^{\Lambda}%
$\textit{ (resp. of }$\mathcal{M}_{(A,\mathcal{E},\mathcal{P})}$%
)\textit{.}\newline$\left(  ii\right)  $~The factor $\mathcal{M}%
_{(A,\mathcal{E},\mathcal{P})}/\mathcal{N}_{(I_{A},\mathcal{E},\mathcal{P})}$
is a presheaf of algebras over the factor ring $\mathcal{C}=A/I_{A}$, with
localization principle, called \emph{presheaf of }$(\mathcal{C},\mathcal{E}%
,\mathcal{P})$\emph{-algebras}\textit{\textit{.}}
\end{propdef}

Remark that, with $(2)$, the constant sheaf $\mathcal{M}_{(A,\mathbb{K}%
,\left\vert .\right\vert )}/\mathcal{N}_{(I_{A},\mathbb{K},\left\vert
.\right\vert )}$\ is exactly equal to $\mathcal{C}=A/I_{A}.$

\begin{notation}
\textit{We denote by }$\left[  \left(  u_{\lambda}\right)  _{\lambda}\right]
_{\mathcal{A}}=\left[  u_{\lambda}\right]  _{\mathcal{A}}$\textit{\ or
}$\left[  u_{\lambda}\right]  $, when no confusion may arise, \textit{the
class of }$\left(  u_{\lambda}\right)  _{\lambda\in\Lambda}$ in\textit{
}$\mathcal{A}(\Omega).$
\end{notation}

\begin{remark}
\label{SPSRemMorphism}We suppose in addition that $\left\{  \left(
a_{\lambda}\right)  _{\lambda}\in A\,\left\vert \,\lim\nolimits_{\Lambda
}a_{\lambda}=0\right.  \right\}  \neq\emptyset$ and that $I_{A}$ satisfies
\begin{equation}
I_{A}\subset\left\{  \left(  a_{\lambda}\right)  _{\lambda}\in A\,\left\vert
\,\lim\nolimits_{\Lambda}a_{\lambda}=0\right.  \right\}  ,
\label{SPSCondAnneau}%
\end{equation}
Then there exists a canonical sheaf embedding of $\mathcal{E}$ into
$\mathcal{A}$ through the morphism of algebra%
\[
\sigma_{\Omega}:\mathcal{E}\left(  \Omega\right)  \rightarrow\mathcal{A}%
\left(  \Omega\right)  ,\ \ f\mapsto\left[  \left(  f\right)  _{\lambda
}\right]  .
\]
Indeed, if $\left[  \left(  f\right)  _{\lambda}\right]  =0$, we have:
$\forall p\in\mathcal{P}(\Omega),\ \left(  p(f)\right)  _{\lambda}%
\in\left\vert I_{A}\right\vert $.\ From (\ref{SPSCondAnneau}), it follows that
$\forall p\in\mathcal{P}(\Omega),\ p(f)=0$.\ Thus $f=0$.
\end{remark}

\begin{remark}
\label{CEP08RemTopoE}For the above algebraic considerations of this section
(and specially Proposition-Definition \ref{SPSPropSubAlg}), we don't need
$\Lambda$ to be a directed set. However, the previous remark shows the
importance of this assumption in order to get non trivial extensions. 
\end{remark}

\section{Overgenerated rings\label{CEP08SecOverTop}}

In almost all works using $(\mathcal{C},\mathcal{E},\mathcal{P})$-type
algebras (see, for examples, \cite{Devoue1,JAM1,JAM2}), the ring $A$ and the
ideal $I_{A}$ are constructed as polynomially overgenerated rings and satisfy
the assumption (\ref{SPSCondAnneau}) of Remark \ref{SPSRemMorphism}. For
$\left(  a_{\lambda}\right)  _{\lambda},\left(  b_{\lambda}\right)  _{\lambda
}\in\mathbb{R}^{\Lambda}$, we shall use the following notation%
\[
a_{\lambda}\ll b_{\lambda}\Leftrightarrow\exists\lambda_{0}\in\Lambda
,\ \forall\lambda\preceq\lambda_{0}:a_{\lambda}\leq b_{\lambda}.
\]

We first give an improved definition of the overgeneration.

\begin{propdef}
\label{Overgenerated}(\textbf{Polynomially overgenerated rings)} Consider
$\mathcal{B}_{0}$ a family of$\ $nets in $(\mathbb{R}_{+}^{\ast})^{\Lambda}$
and $\mathcal{B}$ the subset of elements in $(\mathbb{R}_{+}^{\ast})^{\Lambda
}$ obtained as rational functions with coefficients in $\mathbb{R}_{+}^{\ast}$
of elements in $\mathcal{B}_{0}$ as variables. Set
\[
A_{\mathcal{B}}=\left\{  \left(  a_{\lambda}\right)  _{\lambda}\in
\mathbb{K}^{\Lambda}\mid\exists\left(  b_{\lambda}\right)  _{\lambda}%
\in\mathcal{B}:\left\vert a_{\lambda}\right\vert \ll b_{\lambda}\right\}  .
\]
The set $A_{\mathcal{B}}$ is a solid subring of $\mathbb{K}^{\Lambda}$, called
the ring \emph{(polynomially) overgenerated }by $\mathcal{B}_{0}$ (or by
$\mathcal{B}$).
\end{propdef}

Usually, the set $\mathcal{B}_{0}$ is finite and given by the problem itself.
(See \cite{Devoue1,JAM2}.) The term polynomially refers to the fact that the
growth of elements of $A_{\mathcal{B}}$ is at most polynomial with respect to
the elements of $\mathcal{B}_{0}$. This polynomial overgeneration is
sufficient for the non linearities considered in the quoted references, but,
for example, does not permit to obtain $(\mathcal{C},\mathcal{E},\mathcal{P}%
)$-algebras stable by exponential.

\begin{remark}
With this definition $\mathcal{B}$ is stable by inverse. In many practical
cases and, for example, in the case of Colombeau simplified algebras, which
are a particular case of $(\mathcal{C},\mathcal{E},\mathcal{P})$-algebras,
$\mathcal{B}$ is exactly the set of invertible elements of the ring of
generalized constants.
\end{remark}

As a \textquotedblleft canonical\textquotedblright\ ideal of $A_{\mathcal{B}}%
$, one usually choose%
\begin{equation}
I_{\mathcal{B}}=\left\{  \left(  a_{\lambda}\right)  _{\lambda}\in
\mathbb{K}^{\Lambda}\mid\forall\left(  b_{\lambda}\right)  _{\lambda}%
\in\mathcal{B}:\left\vert a_{\lambda}\right\vert \ll b_{\lambda}\right\}
\text{.} \label{CEP08CanonicIdeal}%
\end{equation}
A routine checking shows that $I_{\mathcal{B}}$ is a solid ideal of
$A_{\mathcal{B}}$. We shall always assume the existence of $\left(
r_{\lambda}\right)  _{\lambda}\in\mathcal{B}$ such that $\lim_{\Lambda
}r_{\lambda}=0$, in order to have (\ref{SPSCondAnneau}) and, thus, the
canonical embedding of $\mathcal{E}\left(  \Omega\right)  $ into
$\mathcal{A}\left(  \Omega\right)  $.\ (This assumption is satisfied in all
practical applications.) We denote by $\mathcal{C}_{\mathcal{B}}%
=A_{\mathcal{B}}/I_{\mathcal{B}}$ the corresponding ring of generalized numbers.

\section{An Austrian Lemma in $(\mathcal{C},\mathcal{E},\mathcal{P}%
)$-algebras\label{SecAL}}

We take here $\mathcal{E}=\mathrm{C}^{\infty}$ with $X=\mathbb{R}^{d}$,
$\mathcal{P}\left(  \mathbb{R}^{d}\right)  $ being the usual family of
seminorms $\left(  P_{K,l}\right)  _{K,l\text{ }}$defined by
\[
P_{K,l}(u)=\sup_{\left\vert \alpha\right\vert \leq l}P_{K,\alpha}(u)\text{
\ with }P_{K,\alpha}(u)=\sup_{x\in K}\left\vert \partial^{\alpha
}u(x)\right\vert \text{, \ }K\subset\subset\Omega,\ \ l=0\text{ or }l=1.
\]
and $\partial^{\alpha}=\dfrac{\partial^{\alpha_{1}+...+\alpha_{d}}}{\partial
z_{1}^{\alpha_{1}}...\partial z_{d}^{\alpha_{d}}}$ for $z=(z_{1},\ldots
,z_{d})\in\Omega$, $l\in\mathbb{N}$, $\alpha=(\alpha_{1},...,\alpha_{d}%
)\in\mathbb{N}^{d}$. We consider a ring of generalized constants
$\mathcal{C}=A_{\mathcal{B}}/I_{\mathcal{B}}$ overgenerated as stated in
Proposition-Definition \ref{Overgenerated}. The ideal $I_{\mathcal{B}}$ is
defined by (\ref{CEP08CanonicIdeal}) and the set of indices $\Lambda$ is
assumed to be left filtering. Recall that%
\begin{align*}
\mathcal{M}(\mathbb{R}^{d}) &  =\mathcal{M}_{(A_{\mathcal{B}},\mathrm{C}%
^{\infty},\mathcal{P})}(\mathbb{R}^{d})=\{\left(  u_{\lambda}\right)
_{\lambda}\in\mathrm{C}^{\infty}(\Omega)^{\Lambda}:\forall p\in\mathcal{P}%
(\mathbb{R}^{d})\text{, }\left(  p(u_{\lambda})\right)  _{\lambda}%
\in\left\vert A_{\mathcal{B}}\right\vert \}\text{,}\\
\mathcal{N}(\mathbb{R}^{d}) &  =\mathcal{N}_{(I_{\mathcal{B}},\mathrm{C}%
^{\infty},\mathcal{P})}(\mathbb{R}^{d})=\{\left(  u_{\lambda}\right)
_{\lambda}\in\mathrm{C}^{\infty}(\Omega)^{\Lambda}:\forall p\in\mathcal{P}%
(\mathbb{R}^{d})\text{, }\left(  p(u_{\lambda})\right)  _{\lambda}%
\in\left\vert I_{\mathcal{B}}\right\vert \}\text{.}%
\end{align*}

\begin{proposition}
\label{CEP080-estimate}Assume that there exists\ $\left(  a_{\lambda}\right)
_{\lambda}\in\mathcal{B}$ with\ $\lim\nolimits_{\Lambda}a_{\lambda}%
=0.\ $Then\textit{ }$\left(  u_{\lambda}\right)  _{\lambda}\in\mathcal{M}%
_{(A_{\mathcal{B}},\mathrm{C}^{\infty},\mathcal{P})}(\mathbb{R}^{d})$ is in
$\mathcal{N}_{(I_{\mathcal{B}},\mathrm{C}^{\infty},\mathcal{P})}%
(\mathbb{R}^{d})$ if, and only if,
\[
\forall K\subset\subset\mathbb{R}^{2},\ \ P_{K,0}\left(  u_{\lambda}\right)
\in\left\vert I_{\mathcal{B}}\right\vert .
\]

\end{proposition}

\begin{remark}
We recall that the set $\mathcal{B}$ is stable by inverse, which could be
assumed for all the $(\mathcal{C},\mathcal{E},\mathcal{P})$-algebras
considered up to now in the literature. Notice also that one has the existence
of $\left(  a_{\lambda}\right)  _{\lambda}\in\mathcal{B}$ such that
$\lim\nolimits_{\Lambda}a_{\lambda}=0$ in all practical cases.
\end{remark}

\begin{proof}
Take $K\subset\subset\Omega$.\ We have to prove that $\forall l\in
\mathbb{N},\;P_{K,l}\left(  u_{\lambda}\right)  \in\left\vert I_{\mathcal{B}%
}\right\vert $.\ By induction, it suffices to prove that $P_{K,0}\left(
u_{\lambda}\right)  \in\left\vert I_{\mathcal{B}}\right\vert $ implies
$P_{K,1}\left(  u_{\lambda}\right)  \in\left\vert I_{\mathcal{B}}\right\vert
$.\ In fact, this amounts to show that $P_{K,0}\left(  u_{\lambda}\right)
\in\left\vert I_{\mathcal{B}}\right\vert $ implies $P_{K,0}\left(  \left(
\partial/\partial x_{i}\right)  u_{\lambda}\right)  \in\left\vert
I_{\mathcal{B}}\right\vert $ for $i\in\left\{  1,\ldots,d\right\}  $. Set
$\delta=\min(1,\operatorname*{dist}(K,\partial\Omega))$ and $L=K+\overline
{B}(0,\delta/2)$. We have $K\subset\subset L\subset\subset\Omega$. Since
$\left(  u_{\lambda}\right)  _{\lambda}\in\mathcal{M}(\mathbb{R}^{d})$, there
exists $\left(  \beta_{\lambda}\right)  _{\lambda}\in\mathcal{B}$ such that
\[
\exists\lambda_{0}\in\Lambda,\forall\lambda\preceq\lambda_{0},\ \ P_{L,2}%
\left(  u_{\lambda}\right)  \leq\beta_{\lambda}.
\]
We may assume that $\lim_{\Lambda}\beta_{\lambda}=+\infty$.\ Indeed, for any
$\left(  \beta_{\lambda}\right)  _{\lambda}\in\mathcal{B}$, we set
$\beta_{\lambda}^{\prime}=a_{\lambda}^{-1}+\beta_{\lambda}$ where $\left(
a_{\lambda}\right)  _{\lambda}\in\mathcal{B}$ is such that $\lim_{\Lambda
}a_{\lambda}=0$.\ Thus, $\lim_{\Lambda}\max\beta_{\lambda}^{\prime}=+\infty$.
Take any $\left(  c_{\lambda}\right)  _{\lambda}\in\mathcal{B}$ and define
$b_{\lambda}=a_{\lambda}c_{\lambda}/\left(  a_{\lambda}+c_{\lambda}\right)  $.
Clearly we have $b_{\lambda}\in\left\vert A_{\mathcal{B}}\right\vert $,
$b_{\lambda}\leq c_{\lambda}$ and $b_{\lambda}\leq a_{\lambda}$.\ Thus
$\lim_{\Lambda}b_{\lambda}=0$.\ Let $\left(  e_{i}\right)  _{1\leq i\leq d}$
be the canonical base of $\mathbb{R}^{d}$. There exists $\lambda_{1}$ such
that, for all $x\in K$, $x+b_{\lambda}\beta_{\lambda}^{-1}e_{i}\in L$ when
$\lambda\preceq\lambda_{1}$, since $\lim_{\Lambda}\beta_{\lambda}^{-1}=0$. By
the Taylor theorem we have, for $x\in K$,%
\[
u_{\lambda}\left(  x+b_{\lambda}\beta_{\lambda}^{-1}e_{i}\right)  =u_{\lambda
}\left(  x\right)  +b_{\lambda}\beta_{\lambda}^{-1}\frac{\partial}{\partial
x_{i}}u_{\lambda}\left(  x\right)  +\frac{1}{2}\left(  b_{\lambda}%
\beta_{\lambda}^{-1}\right)  ^{2}\frac{\partial^{2}}{\partial x_{i}^{2}%
}u_{\lambda}\left(  x+\theta b_{\lambda}\beta_{\lambda}^{-1}e_{i}\right)
\]
with $0\leq\theta\leq1$. It follows that%
\[
\frac{\partial}{\partial x_{i}}u_{\lambda}\left(  x\right)  =b_{\lambda}%
^{-1}\beta_{\lambda}\left(  u_{\lambda}\left(  x+b_{\lambda}\beta_{\lambda
}^{-1}e_{i}\right)  -u_{\lambda}\left(  x\right)  \right)  -\frac{1}{2}\left(
b_{\lambda}\beta_{\lambda}^{-1}\right)  \frac{\partial^{2}}{\partial x_{i}%
^{2}}u_{\lambda}\left(  x+\theta b_{\lambda}\beta_{\lambda}^{-1}e_{i}\right)
.
\]
Thus%
\[
\left\vert \frac{\partial}{\partial x_{i}}u_{\lambda}\left(  x\right)
\right\vert \leq2b_{\lambda}^{-1}\beta_{\lambda}P_{L,0}\left(  u_{\lambda
}\right)  +\frac{1}{2}b_{\lambda}\beta_{\lambda}^{-1}P_{L,2}\left(
u_{\lambda}\right)  \leq2b_{\lambda}^{-1}\beta_{\lambda}P_{L,0}\left(
u_{\lambda}\right)  +\frac{1}{2}b_{\lambda}%
\]
for $\lambda\preceq\lambda_{2}$ with $\lambda_{2}\preceq\lambda_{j}$, $0\leq
j\leq1$. As $P_{K,0}\left(  u_{\lambda}\right)  \in\left\vert I_{\mathcal{B}%
}\right\vert $, we have $P_{L,0}\left(  u_{\lambda}\right)  \leq
(1/4)b_{\lambda}^{2}\beta_{\lambda}^{-1}\in\mathcal{B}$ for $\lambda
\preceq\lambda_{3}$ for some $\lambda_{3}$. Thus\vspace{-0.04in}
\[
\left\vert \frac{\partial}{\partial x_{i}}u_{\lambda}\left(  x\right)
\right\vert \leq b_{\lambda}\ \ \text{for }\lambda\preceq\lambda_{4}\text{
with }\lambda_{4}\preceq\lambda_{j}\text{,\ }3\leq j\leq4.
\]
Finally, $P_{K,0}\left(  \left(  \partial/\partial x_{i}\right)  u_{\lambda
}\right)  \in\left\vert I_{\mathcal{B}}\right\vert $ as expected.
\end{proof}

\end{document}